\documentclass[12pt]{article}
\usepackage{times,eepic}
\makeatletter
\def\address{\m@th\@ifnextchar[\@address{\@address[]}}

\def\@address[#1]#2{
\expandafter\def\expandafter\@addressname\expandafter
{\@addressname{
  \adr{#1}\ \parbox[t]{5in}{
     \ignorespaces #2}\par }}}
\def\@addressname{}
\def\adr#1{{\normalsize\unskip$^{#1}$}}
\def\@maketitle{%
\def\and{{\rm and}}
  \newpage
  \null
  {\centering
  \let \footnote \thanks
    {\LARGE\bf   \@title \par}%
    \vskip 1.5em%
      \lineskip .5em%
    {\sc\large   \@author\par}
      \vspace{1em}
    {\small \@addressname}

  }%
  \par
  \vskip 1.5em}
\makeatother
\title{Quantum Moduli Spaces of Flat Connections}
\author{Anton Yu. Alekseev \adr{1} 
        \and \ Volker Schomerus \adr{2}}  
\address[1]{Institute of Theoretical Physics, Uppsala University,  
          Box 803 S-75108, \\ Uppsala, Sweden; \ 
          e-mail: alekseev@teorfys.uu.se }   
\address[2]{II. Institut f\"ur Theoretische Physik, Universit\"at Hamburg,
          Luruper Chaussee 149,\\ 22761 Hamburg, Germany; \  
          e-mail: vschomer@x4u2.desy.de }
\begingroup\makeatletter\ifx\SetFigFont\undefined
\def\x#1#2#3#4#5#6#7\relax{\def\x{#1#2#3#4#5#6}}%
\expandafter\x\fmtname xxxxxx\relax \def\y{splain}%
\ifx\x\y   % LaTeX or SliTeX?
\gdef\SetFigFont#1#2#3{%
  \ifnum #1<17\tiny\else \ifnum #1<20\small\else
  \ifnum #1<24\normalsize\else \ifnum #1<29\large\else
  \ifnum #1<34\Large\else \ifnum #1<41\LARGE\else
     \huge\fi\fi\fi\fi\fi\fi
  \csname #3\endcsname}%
\else
\gdef\SetFigFont#1#2#3{\begingroup
  \count@#1\relax \ifnum 25<\count@\count@25\fi
  \def\x{\endgroup\@setsize\SetFigFont{#2pt}}%
  \expandafter\x
    \csname \romannumeral\the\count@ pt\expandafter\endcsname
    \csname @\romannumeral\the\count@ pt\endcsname
  \csname #3\endcsname}%
\fi
\fi\endgroup
\renewcommand{\sloppy}{\tolerance=5000}
\newcommand{\sodot}{\,{\scriptscriptstyle \odot}\, } 
\def\sg{{\sf g}}
\makeatletter
\renewcommand\section{\@startsection{section}{1}{\z@}%
                                     {-3.25ex\@plus -1ex \@minus -.2ex}%
                                     {1.5ex \@plus .2ex}%
                                     {\normalfont\large\bfseries}}
\renewcommand\subsection{\@startsection{subsection}{2}{\z@}%
                                     {-3.25ex\@plus -1ex \@minus -.2ex}%
                                     {1.5ex \@plus .2ex}%
                                     {\normalfont\normalsize\bfseries}}
\makeatother
\newtheorem{theo}{Theorem}

\newtheorem{prop}[theo]{Proposition}  
  
\def\sU{{\sf U}}        \def\sT{{\sf T}}     \def\sG{{\sf G}}
\def\sM{{\sf M}}        \def\sv{{\sf v}}     \def\sA{{\sf A}}
\def\sB{{\sf B}}
\def\rR{{\rm R}}
     
     \def\tr{{\mbox{\it tr\/}}}
\def\be{\begin{equation}}      \def\ba{\begin{eqnarray}}  
\def\ee{\end{equation}}        \def\ea{\end{eqnarray}}  
\def\o{\otimes }  
  
\def\D{\Delta }  
  
\def\cH{{\cal H}}  

\def\cF{{\cal F}}
\def\cJ{{\cal J}}
\def\cC{{\cal C}}
  
\newcommand{\cL}{{\cal L}}  
  
\def\cA{{\cal A}}  
  
\def\cG{{\cal G}}  
\def\F{{\cal F}}
  
\def\ti{\times }  
  
\def\s{\sigma }  
\def\sg{{\rm g}}

\def\a{\alpha }          \def\b{\beta }           
                   
          \def\t{\tau}  
\def\nn{\nonumber}  
\def\sss{\scriptscriptstyle}
\newcommand{\Ae}[1]{\stackrel{\sss #1}{\raisebox{.1pt}{\sf A}}}
\newcommand{\Be}[1]{\stackrel{\sss #1}{\raisebox{.1pt}{\sf B}}}
\newcommand{\U}[1]{\stackrel{\sss #1}{\raisebox{.1pt}{\sf U}}}
\newcommand{\M}[1]{\stackrel{\sss #1}{\raisebox{.1pt}{\sf M}}}
\newcommand{\T}[1]{\stackrel{\sss #1}{\raisebox{.1pt}{\sf T}}}
\newcommand{\G}[1]{\stackrel{\sss #1}{\raisebox{.1pt}{\sf G}}}
\begin{document}  

\maketitle

\begin{abstract} 
Using the formalism of discrete quantum group gauge theory, 
one can construct the quantum algebras of observables for the 
Hamiltonian Chern-Simons model. The resulting {\em moduli algebras} 
provide quantizations of the algebra of functions on the moduli spaces
of flat connections on a punctured 2-dimensional surface.  In this note we 
describe some features of these moduli algebras with special 
emphasis on the natural action of mapping class groups. This leads, 
in particular, to a closed formula for representations of the 
mapping class groups on conformal blocks.    
\end{abstract}  \vspace*{-13cm}
{\tt \hfill DESY 96-261}\\
{\tt \phantom{a} \hfill q-alg/9612037}\\[12cm]
\setcounter{footnote}{1} 
\footnotetext{To be published in the Proceedings of {\em XXI International 
Colloquium on Group Theoretical Methods in Physics, Symposium on 
Quantum Groups}, Goslar 1996, edited by H.-D. Doebner and V.K. Dobrev, 
Heron Press (Sofia).} 

\section{Introduction}   

Moduli spaces of flat connections on a Riemann surface $\Sigma$
enter physics as phase spaces of the so-called Chern-Simons gauge
theory on a three dimensional manifold $\Sigma \ti {\bf R }$. They
possess a famous Poisson structure that has been studied extensively
in the mathematical literature (see e.g. \cite{PB}). The information
about a particular flat connection on $\Sigma$ can be encoded in its
{\em holonomies along nontrivial closed curves}. Hence, one may
regard points on the moduli spaces as classes of homomorphisms
{}from the fundamental group $\pi_1(\Sigma)$ into the gauge group
$G$. Such a description based on $\pi_1(\Sigma)$ furnishes a
natural action of the {\em mapping class group} of $\Sigma$ on
moduli spaces. The latter is known to respect the canonical
Poisson structure. 

Different methods have been employed to quantize the Hamiltonian
Chern Simons theory (see e.g. \cite{geo,AGS}). The combinatorial
quantization developed in \cite{AGS} resulted in a noncommutative
algebra $\cA_{\sf CS}$ of `functions on the quantum moduli space
of flat connections'. It is one feature of this approach that it
follows closely the classical construction of moduli spaces with
certain {\em quantum holonomies} replacing their classical
counterparts (see below). The deformed universal enveloping algebras
$\cG = U_q(\sg)$ play the role of gauge symmetries in the quantized
theory. As it is expected, the action of the mapping class group
survives quantization and gives rise to automorphisms of the
algebras $\cA_{\sf CS}$ of observables in the quantum theory. 

In this note we reconstruct quantum moduli spaces of flat
connections from the described properties, namely that (1) they
are built up from quantum holonomies and (2) they admit an action
of the mapping class group. The interplay of classical and quantum
symmetries (i.e. the action of mapping class groups and deformed
universal enveloping algebras) naturally leads to the algebras
$\cA_{\sf CS}$ which were found in \cite{AGS,AlSc} by quantizing
the Poisson structure on moduli spaces of flat connections. There
the analysis was essentially based on a lattice formulation of the
Hamiltonian Chern-Simons theory due to V.V. Fock and A.A. Rosly
\cite{PB}. Even without making this connection to the classical
Chern-Simons theory explicit, our constructions below have an
immediate mathematical application: they provide a closed and very
elegant formula for the Reshetikhin-Turaev representations of
mapping class groups.

\section{Products on the Dual of $U_q(\sg)$}
 
The space $\cG'$ of linear forms on the quantum universal enveloping
algebras $\cG = U_q(\sg)$ can be equipped with a number of different
multiplications. In this section we describe three such product
structures and discuss some of their properties thereby discovering
natural quantum analogues of holonomies along open and closed curves. 

\subsection{Quantum group algebras} It is well known that the 
Hopf algebra structure of $\cG$ admits to define a natural product 
$a \circ b$ of two linear forms $a,b: \cG \rightarrow \cC$. By 
definition, we evaluate $a \circ b$ on elements $\xi \in \cG$ 
with the help of  $(a \circ b)(\xi) = (a \o b) (\D(\xi))$. The 
associative algebra $\cF$ which emerges from this simple construction 
is known as {\em quantum group algebra}. Representations $\tau$ of 
$\cG$ may be reinterpreted as matrices of linear forms on $\cG$ 
and hence -- with the natural product on $\cG'$ -- as matrices 
with elements in the non-commutative quantum group algebra $\cF$. 
Let us consider one particular representation $\t$ and denote
the associated $\cF$-valued matrix by $\sU$. From the definition
of $\circ$ in terms of $\D$ and the famous intertwining property
$\rR \D(.) = \D'(.) \rR$ one deduces the following quadratic
relations:  
\be  \label{U}     \rR  \ \U{1}{}\  \U{2}{} \ = \ \U{2}{} \
         \U{1}{}\  \rR \ \ .
\ee
Here $\U{1}{}  = \sU \o {\bf 1}, \U{2}{} = {\bf 1} \o \sU$
and multiplication of matrix elements with $\circ$ is understood.  
$\rR$ denotes the universal $\rR$-matrix of $\cG$ evaluated in the 
representation $\t$. In the standard example of $\cG = U_q({\sf sl}_2)$
we may fix $\t$ to be the two dimensional fundamental representation.
Then the quadratic relations (\ref{U}) become equations between
products of algebra valued $4 \ti 4$ matrices and $\rR$ is the
famous $4 \ti 4$-matrix-solution of the Yang-Baxter equation.
When, in addition to these quadratic relations, the value of the
{\em quantum determinant} is fixed to $1$, we obtain a complete
description of $\cF$. We disregard such determinant relations
in the following to simplify our discussion. For a precise
treatment it is preferable to assemble the matrices $\sU$ into
{\em universal elements}. The latter obey certain {\em
functoriality relations} which imply quadratic relations and
constrain the quantum determinant at the same time. This and
other topics, such as the existence of a $\ast$-operation,
are explained in \cite{AGS,AlSc}. 

To investigate the quantum symmetries of $\F$ we introduce another
copy $\sT$ of our matrix $\sU$ such that 
\be \label{T} \rR  \ \T{1}{} \ \T{2}{} \ =\ 
    \ \T{2}{}\  \T{1}{}\  \rR \ \ \ee
and matrix elements of $\sT$ commute with matrix elements of
$\sU$. Together, matrix elements of $\sU, \sT$ generate the
algebra $\cF \o \cF$ and we can use $\sT$ in defining  
a homomorphism  $\Phi: \F \rightarrow \F \o \F$,   
$$     \Phi : \ \sU \ \mapsto \ \sU\  \sT \ \ . $$
On the right hand side, $\sU\sT$ is an ordinary product of
algebra valued matrices. Even though it is 
quite elementary, let us demonstrate once that $\Phi$ indeed
extends to a homomorphism. Similar computations will then be
left to the reader in what follows.
\ba 
 \Phi (\ \rR \ \U{1}{}\U{2}{})  & = & 
      \rR \ \U{1}{}\T{1}{} \U{2}{} \T{2}{} 
     \ = \ \rR \ \U{1}{} \U{2}{} \T{1}{}\T{2}{} 
     \  = \  \ \U{2}{}\U{1}{}\T{2}{}\T{1}{} \ \rR \nn\\[1mm] 
      & = & \U{2}{} \T{2}{} \U{1}{} \T{1}{}\  \rR 
     \ = \ \Phi (\ \U{2}{} \U{1}{}\  \rR ) \ . \nn 
\ea
We call the above symmetry $\Phi = \Phi_{\sf R}$ a {\em right
regular transformation} because objects $\sT$ are multiplied from
the right. 
\subsection{Link- and loop-algebras}
One would certainly expect to find a similar left regular transformation
$\Phi_{\sf L}$ which maps $\sU$ to $\sT^{-1}_{\sf L} \sU$. A short
computation reveals that such a map gives rise to a homomorphism,
if we replace $\rR$ by $\rR' = P \rR P$ in the exchange relations
(\ref{T}) for $\sT= \sT_{\sf L}$. Let us follow a slightly different
route here: we will leave the relations for $\sT$ untouched and
manipulate the exchange relations for $\sU$ instead so that the
resulting algebra has the desired two commuting quantum symmetries
(QS) under left/right-regular transformations with $\sT$. To be more
concrete, we introduce the {\em link-algebra} $\cJ$ that is
generated by elements of an algebra valued matrix $\sG$ such that\\ 
\parbox[b]{25em}{
\ba     \hspace*{-3em}
        \rR' \ \G{1}{} \ \G{2}{} \ &=& \ \G{2}{} \ \G{1}{}\ \ \rR \
        \label{G} \\[2mm]
        \cJ \ \mbox{has QS: }\ \ \Phi_{\sf R}:\  \sG   
        & & \hspace*{-1.5em} \mapsto \ \sG \ \sT_{\sf R}\ \ , \  
        \ \Phi_{\sf L}: \sG \ \mapsto \ \sT^{-1}_{\sf L} \sG \ \ . \nn
\ea } \hspace*{0.5em}
\raisebox{4ex}{{
\setlength{\unitlength}{0.00038333in}
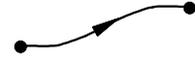
\begin{picture}(2480,709)(0,-10)
\thicklines
\path(54,100)   (108.832,97.490)
        (159.833,95.314)
        (207.196,93.474)
        (251.115,91.968)
        (291.783,90.796)
        (329.395,89.959)
        (396.222,89.290)
        (453.147,89.959)
        (501.719,91.968)
        (543.487,95.314)
        (580.000,100.000)

\path(580,100)  (620.802,108.794)
        (669.523,122.260)
        (723.410,139.165)
        (779.711,158.271)
        (835.673,178.344)
        (888.544,198.146)
        (935.570,216.444)
        (974.000,232.000)

\path(974,232)  (1032.416,258.929)
        (1105.430,295.715)
        (1177.979,333.393)
        (1235.000,363.000)

\path(1235,363) (1292.974,392.315)
        (1366.050,429.793)
        (1439.351,466.373)
        (1498.000,493.000)

\path(1498,493) (1536.312,508.724)
        (1583.254,527.140)
        (1636.074,547.015)
        (1692.019,567.111)
        (1748.337,586.194)
        (1802.275,603.026)
        (1851.080,616.374)
        (1892.000,625.000)

\path(1892,625) (1928.398,629.853)
        (1970.083,633.319)
        (2018.604,635.399)
        (2075.511,636.092)
        (2142.355,635.399)
        (2179.987,634.533)
        (2220.684,633.319)
        (2264.640,631.759)
        (2312.049,629.853)
        (2363.105,627.600)
        (2418.000,625.000)

\blacken\path(1384,446)(1081,347)(1137,244)(1384,446)
\path(1384,446)(1081,347)(1137,244)(1384,446)
\put(2400,615){\blacken\ellipse{144}{144}}
\put(2400,615){\ellipse{144}{144}}
\put(80,79){\blacken\ellipse{144}{144}}
\put(80,79){\ellipse{144}{144}}
\end{picture} }}\\[-1mm] 
Here $\sT_{\sf R}, \sT_{\sf L}$ obey the exchange relations (\ref{T})
and commute with each other and with matrix elements of $\sG$. The
difference between eqs. (\ref{U}) and (\ref{G}) is that one of
the $\rR$-matrices has been replaced by $\rR'$. We wish to emphasize
that the precise meaning of eq. (\ref{G}) is the same as in eq.
(\ref{U}). Namely, the think of $\cJ$ as being generated by elements
of $\cG'$ subject to a new multiplication $(a \sodot b)(\xi) :=
(a \o b)(\rR \D(\xi))$. One may check that this prescription defines
an associative product on $\cG'$ such that matrix elements of the
representations obey the exchange relations (\ref{G}). 

The algebra $\cJ$ that is generated by the components of $\sG$ was 
constructed so that it admits two commuting quantum symmetries
which we called left- and right-regular transformations. Clearly,
these symmetries are remnants of the usual left- and right regular
actions of a group on itself. Classical groups, however, possess
another important symmetry: the adjoint action, i.e. the action
of the group on itself by conjugation. It is easy to see that
this symmetry is lost when we pass from the group to the quantum
object $\sG$. In other words, the map $\sG \mapsto \sT^{-1} \sG
\sT $ does not extend to a homomorphism of algebras, if $\sT$ and
$\sG$ obey the relations we have specified above. But from our
experience with $\sG$ we know already, how a manipulation
of the fundamental exchange relations may influence the quantum
symmetry of the corresponding noncommutative algebra. It is not
too hard to come up with  new relations for an algebra valued
matrix $\sM$ which are consistent with the adjoint transformation
that we were looking for, \\
\parbox[b]{25em}{
\ba    \label{M}
    \rR'\ \M{1}{}\  \rR \ \M{2}{} \ &=& \ \M{2}{}\ \rR'\ \M{1}{}
   \ \rR \ \  \\[2mm]
\cL \ \mbox{ has QS: } \ & &  \Phi_{\sf Ad} :\  \sM \ \mapsto\
     \sT^{-1} \sM\  \sT \nn\ \ .
\ea } \hspace*{1em}
\raisebox{4ex}{
\setlength{\unitlength}{0.00043333in}
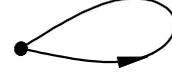
\begin{picture}(2002,900)(0,-10)
\thicklines
\path(80,236)   (121.175,226.525)
        (161.650,217.328)
        (201.433,208.407)
        (240.535,199.762)
        (278.963,191.394)
        (316.728,183.300)
        (390.304,167.939)
        (461.335,153.676)
        (529.895,140.507)
        (596.057,128.431)
        (659.896,117.443)
        (721.484,107.542)
        (780.896,98.725)
        (838.204,90.988)
        (893.482,84.328)
        (946.804,78.743)
        (998.243,74.230)
        (1047.873,70.787)
        (1095.767,68.409)
        (1141.999,67.094)
        (1186.643,66.840)
        (1229.771,67.643)
        (1271.457,69.501)
        (1311.775,72.410)
        (1350.798,76.369)
        (1425.255,87.420)
        (1495.416,102.632)
        (1561.868,121.981)
        (1625.200,145.445)
        (1686.000,173.000)

\path(1686,173) (1735.854,200.550)
        (1787.683,235.451)
        (1838.323,277.146)
        (1884.613,325.076)
        (1923.388,378.684)
        (1951.486,437.413)
        (1965.744,500.704)
        (1963.000,568.000)

\path(1963,568) (1942.360,632.059)
        (1907.137,686.484)
        (1860.499,731.853)
        (1805.616,768.740)
        (1745.656,797.722)
        (1683.788,819.376)
        (1623.179,834.276)
        (1567.000,843.000)

\path(1567,843) (1500.585,848.010)
        (1433.174,848.331)
        (1364.224,843.741)
        (1293.190,834.018)
        (1219.527,818.938)
        (1181.540,809.320)
        (1142.691,798.280)
        (1102.913,785.790)
        (1062.137,771.822)
        (1020.296,756.348)
        (977.321,739.341)
        (933.145,720.773)
        (887.698,700.615)
        (840.914,678.841)
        (792.724,655.422)
        (743.059,630.331)
        (691.853,603.539)
        (639.036,575.020)
        (584.542,544.744)
        (528.301,512.685)
        (470.245,478.815)
        (410.308,443.106)
        (348.419,405.530)
        (284.513,366.059)
        (218.519,324.665)
        (150.371,281.322)
        (80.000,236.000)

\blacken\path(1571,116)(1272,120)(1291,12)(1571,116)
\path(1571,116)(1272,120)(1291,12)(1571,116)
\put(76,223){\blacken\ellipse{136}{136}}
\put(76,223){\ellipse{136}{136}}
\put(80,236){\blacken\ellipse{136}{136}}
\put(80,236){\ellipse{136}{136}}
\end{picture}

}\\[-1mm]
Matrix elements of $\sM$ generate the {\em loop-algebra} $\cL$ and 
we call $\Phi_{\sf Ad}$ the {\em adjoint transformation}. Again,
there exits a product on the dual of $\cG$ which implies that
representation matrices obey eq. (\ref{M}) and, strictly speaking,
it is $\cG'$ endowed with this product that we denote by $\cL$. The
relations which we have just constructed are not new. In fact, they
appeared in a work of Reshetikhin and Semenov-Tian-Shansky on the
deformations of universal enveloping algebras. For finite dimensional
semisimple modular Hopf algebras $\cG$, the loop-algebra $\cL$ is
even isomorphic to $\cG$ \cite{AlSc}. 

Elements in the loop-algebra which transform trivially under the 
symmetry $\Phi_{\sf Ad}$ form a subalgebra $\cL^{\cG} \subset 
\cL$. \footnote{In more mathematical terms, $\Phi_{\sf Ad}$ is called 
a {\em co-action} and elements in $\cL^{\cG}$ are said to be {\em 
co-invariant} with respect to $\Phi_{\sf Ad}$}. 
The {\em quantum trace} $\tr_q$ is especially 
designed to produce such an element $c = \tr_q(\sM) \in \cL^\cG$ 
from the algebra valued matrix $\sM$. We will see more examples 
of this construction later when we consider algebras composed 
from several copies of $\cL$. In our particular case here, the 
element $c = \tr_q(\sM)$ is central, i.e. it commutes with all 
the matrix elements of $\sM$. Further computations reveal that 
$c$ generates the {\em fusion algebra} (or {\em Verlinde algebra} 
\cite{Ver}) which is defined by the multiplicities in the Clebsch 
Gordan decomposition of $\cG = U_q(\sg)$.   
 
The behaviour of $\sG$ and $\sM$ under quantum symmetry
transformations may be attributed to different geometric situations.
In some sense, $\sG$ behaves very much like a holonomy along an open
curve ({\em link}) on which gauge transformations at the two endpoints
act independently from the left and from the right. The quantum
symmetry of $\sM$ should be compared with the transformation law of
holonomies along a closed curve ({\em loop}). On such holonomies,
gauge transformations act by conjugation. Our pictorial presentations
display these differences and encode the symmetry properties of $\sG$
and $\sM$.      
                   
\section{Braided Tensor Products and Braid Groups}
 
In the previous  section we met three quite different algebraic 
structures. From now on, we concentrate on the loop-algebra and
use it as a basic building block to construct more complicated 
algebras. New classical symmetries will emerge in addition to  
the quantum symmetries that we have discussed already.  

\subsection{Multiloop-algebras $\cL_N$} Whereas the tensor product  
is one of the most fundamental constructions in algebra, it is well   
known to break quantum symmetries and hence becomes useless in our 
present context. If we want to combine several copies of $\cL$ into 
one algebra while preserving the quantum symmetry, we are forced to 
use {\em braided tensor products} (see e.g. \cite{maj}). The $N$-fold
braided tensor product  $\cL_N$ of the loop-algebra $\cL$ is generated
by matrix elements  of $\sM_\nu, \nu = 1, \dots, N$ satisfying eq. 
(\ref{M}) and the exchange relations\\[-0mm] 
\parbox[b]{22em}{
\ba   \rR^{-1} \M{1}_\nu\  \rR \ \M{2}_\mu  & = &  
     \M{2}_\mu \ \rR^{-1} \M{1}_\nu \  \rR\ \ \  
      \ \  \nu < \mu \nn\\[4mm]
      \cL_N \mbox{ has QS:}& & \Phi_{\sf Ad}:\     
       \sM_\nu \ \mapsto \ \sT^{-1} \sM_\nu \ \sT \ \ . \nn 
\ea} \hspace*{-1em}
\raisebox{1ex}{
\setlength{\unitlength}{0.00043333in}
\begin{picture}(3861,2186)(0,-10)
\thicklines
\put(949,1373){\blacken\ellipse{42}{42}}
\put(949,1373){\ellipse{42}{42}}
\put(1269,1629){\blacken\ellipse{42}{42}}
\put(1269,1629){\ellipse{42}{42}}
\put(1652,1757){\blacken\ellipse{42}{42}}
\put(1652,1757){\ellipse{42}{42}}
\put(3060,1117){\blacken\ellipse{42}{42}}
\put(3060,1117){\ellipse{42}{42}}
\put(3654,892){\makebox(0,0)[lb]{\smash{{{\SetFigFont{12}{14.4}{rm}1}}}}}
\put(204,1192){\makebox(0,0)[lb]{\smash{{{\SetFigFont{12}{14.4}{rm}N}}}}}
\put(2229,2092){\makebox(0,0)[lb]{\smash{{{\SetFigFont{12}{14.4}{rm}{$\nu$}}}}}}
\path(1926,443) (1964.915,426.484)
        (2003.188,410.364)
        (2077.845,379.301)
        (2150.042,349.796)
        (2219.852,321.833)
        (2287.346,295.396)
        (2352.595,270.470)
        (2415.673,247.039)
        (2476.650,225.089)
        (2535.600,204.603)
        (2592.592,185.566)
        (2647.701,167.963)
        (2700.997,151.778)
        (2752.552,136.995)
        (2802.438,123.600)
        (2850.728,111.576)
        (2897.493,100.908)
        (2942.804,91.580)
        (2986.734,83.578)
        (3029.356,76.885)
        (3070.739,71.486)
        (3110.958,67.366)
        (3150.082,64.509)
        (3188.185,62.900)
        (3225.338,62.522)
        (3297.083,65.402)
        (3365.892,73.024)
        (3432.339,85.265)
        (3497.000,102.000)

\path(3497,102) (3551.049,120.500)
        (3608.179,145.917)
        (3665.219,178.239)
        (3718.995,217.454)
        (3766.334,263.549)
        (3804.063,316.514)
        (3829.009,376.335)
        (3838.000,443.000)

\path(3838,443) (3828.665,509.857)
        (3803.419,569.685)
        (3765.444,622.517)
        (3717.921,668.391)
        (3664.034,707.341)
        (3606.963,739.402)
        (3549.891,764.610)
        (3496.000,783.000)

\path(3496,783) (3431.561,799.502)
        (3365.313,811.555)
        (3296.680,819.035)
        (3225.088,821.816)
        (3188.003,821.404)
        (3149.963,819.771)
        (3110.895,816.901)
        (3070.729,812.777)
        (3029.392,807.385)
        (2986.812,800.708)
        (2942.918,792.731)
        (2897.637,783.437)
        (2850.899,772.813)
        (2802.631,760.841)
        (2752.760,747.506)
        (2701.216,732.793)
        (2647.927,716.685)
        (2592.821,699.168)
        (2535.825,680.225)
        (2476.868,659.841)
        (2415.879,637.999)
        (2352.785,614.686)
        (2287.514,589.883)
        (2219.995,563.577)
        (2150.156,535.752)
        (2077.925,506.391)
        (2003.230,475.479)
        (1964.937,459.436)
        (1926.000,443.000)

\put(1920,431){\blacken\ellipse{136}{136}}
\put(1920,431){\ellipse{136}{136}}
\blacken\path(559,1055)(839,949)(859,1057)(559,1055)
\path(559,1055)(839,949)(859,1057)(559,1055)
\path(1920,431) (1849.763,476.904)
        (1781.738,520.799)
        (1715.857,562.715)
        (1652.050,602.679)
        (1590.251,640.720)
        (1530.390,676.865)
        (1472.400,711.144)
        (1416.212,743.584)
        (1361.758,774.214)
        (1308.970,803.061)
        (1257.779,830.155)
        (1208.118,855.523)
        (1159.919,879.194)
        (1113.112,901.195)
        (1067.629,921.556)
        (1023.404,940.304)
        (980.366,957.467)
        (938.449,973.075)
        (897.583,987.154)
        (857.701,999.734)
        (818.734,1010.842)
        (780.615,1020.507)
        (706.644,1035.621)
        (635.243,1045.302)
        (565.866,1049.774)
        (497.967,1049.265)
        (431.000,1044.000)

\path(431,1044) (374.696,1034.574)
        (314.240,1018.798)
        (252.732,996.176)
        (193.276,966.211)
        (138.973,928.408)
        (92.925,882.269)
        (58.233,827.298)
        (38.000,763.000)

\path(38,763)   (35.849,696.309)
        (50.350,633.616)
        (78.395,575.468)
        (116.876,522.410)
        (162.685,474.988)
        (212.714,433.749)
        (263.855,399.238)
        (313.000,372.000)

\path(313,372)  (373.799,344.439)
        (437.134,320.947)
        (503.592,301.545)
        (573.765,286.255)
        (648.240,275.096)
        (687.275,271.073)
        (727.607,268.090)
        (769.310,266.151)
        (812.456,265.258)
        (857.121,265.413)
        (903.376,266.620)
        (951.297,268.880)
        (1000.956,272.197)
        (1052.428,276.572)
        (1105.786,282.010)
        (1161.103,288.511)
        (1218.454,296.079)
        (1277.911,304.716)
        (1339.550,314.425)
        (1403.443,325.209)
        (1469.663,337.070)
        (1538.286,350.011)
        (1609.383,364.033)
        (1683.030,379.141)
        (1720.832,387.103)
        (1759.299,395.336)
        (1798.440,403.843)
        (1838.265,412.622)
        (1878.781,421.674)
        (1920.000,431.000)

\put(1898,448){\blacken\ellipse{136}{136}}
\put(1898,448){\ellipse{136}{136}}
\blacken\path(3374,66)(3080,122)(3080,12)(3374,66)
\path(3374,66)(3080,122)(3080,12)(3374,66)
\put(1926,443){\blacken\ellipse{136}{136}}
\put(1926,443){\ellipse{136}{136}}
\path(1898,448) (1964.936,498.426)
        (2029.240,547.455)
        (2090.962,595.140)
        (2150.151,641.536)
        (2206.857,686.698)
        (2261.130,730.679)
        (2313.018,773.534)
        (2362.572,815.317)
        (2409.840,856.083)
        (2454.872,895.887)
        (2497.718,934.781)
        (2538.427,972.822)
        (2577.049,1010.063)
        (2613.632,1046.559)
        (2680.882,1117.531)
        (2740.574,1186.175)
        (2793.102,1252.924)
        (2838.863,1318.214)
        (2878.252,1382.480)
        (2911.664,1446.156)
        (2939.496,1509.679)
        (2962.143,1573.482)
        (2980.000,1638.000)

\path(2980,1638)        (2990.824,1693.910)
        (2997.182,1755.916)
        (2997.510,1821.270)
        (2990.245,1887.220)
        (2973.825,1951.018)
        (2946.686,2009.914)
        (2907.265,2061.158)
        (2854.000,2102.000)

\path(2854,2102)        (2791.717,2127.277)
        (2727.549,2135.304)
        (2663.070,2128.863)
        (2599.850,2110.735)
        (2539.461,2083.703)
        (2483.476,2050.548)
        (2433.465,2014.053)
        (2391.000,1977.000)

\path(2391,1977)        (2344.354,1929.315)
        (2300.661,1877.858)
        (2259.742,1822.071)
        (2221.416,1761.392)
        (2185.503,1695.264)
        (2151.821,1623.124)
        (2135.760,1584.626)
        (2120.189,1544.415)
        (2105.087,1502.421)
        (2090.429,1458.575)
        (2076.193,1412.806)
        (2062.358,1365.046)
        (2048.899,1315.222)
        (2035.796,1263.266)
        (2023.024,1209.108)
        (2010.562,1152.678)
        (1998.387,1093.905)
        (1986.477,1032.720)
        (1974.808,969.053)
        (1963.358,902.833)
        (1952.105,833.991)
        (1941.026,762.457)
        (1930.099,688.161)
        (1924.685,649.955)
        (1919.301,611.033)
        (1913.943,571.385)
        (1908.608,531.003)
        (1903.295,489.877)
        (1898.000,448.000)

\blacken\path(2949,1514)(2754,1286)(2849,1231)(2949,1514)
\path(2949,1514)(2754,1286)(2849,1231)(2949,1514)
\end{picture} }\\
An ordinary tensor product of loop-algebras would mean to replace
all $\rR$'s by unit-matrices in the exchange relations of $\sM_\nu$
with $\sM_\mu$. The reader is invited to check that such commutation 
relations are inconsistent with the quantum symmetry. On the 
other hand, $\Phi_{\sf Ad}$ preserves the exchange relations with 
nontrivial $\rR$-matrices. We call $\cL_N$ a {\em multiloop-algebra}
and write $\cL_N^{\cG}$ for the algebra of elements in $\cL_N$ which 
transform trivially under $\Phi_{\sf Ad}$. Note
that the picture we have drawn to illustrate the $N$-fold braided 
tensor product is in agreement with our previous rules. In particular, 
it has one vertex keeping track of the quantum symmetry and $N$ loops 
which symbolize the $N$ loop-algebras that generate $\cL_N$.
 
It turns out that the multiloop-algebras $\cL_N$ have interesting 
classical symmetries in addition to the quantum symmetry $\Phi_{
\sf Ad}$. In fact, elements of the braid group $B_N$ on $N$ strands
act as automorphisms on $\cL_N$. Let $\s_\rho$ denote the generator 
which corresponds to an exchange of the $\rho^{th}$ with the $(\rho
+1)^{st}$ strand. Their action on the multiloop-algebras $\cL_N$
may be defined by 
\ba
\s_\rho(\sM_\nu)\  = \ \sM_\nu \ \ & & \mbox{ for all } \ \ \nu \neq 
   \rho,\rho+1  \\[3mm] 
\s_\rho(\sM_\rho) \ = \ \sM_{\rho+1} \ \ \ \ & , &
    \s_\rho(\sM_{\rho+1})
    \ \sim\   \sM_{\rho+1}^{-1} \ \sM_\rho \ \sM_{\rho+1}\ \ . 
\ea 
This is consistent with all exchange relations of $\cL_N$ and  
respects the defining relations in $B_N$. The `$\sim$' in the 
second line means `up to a scalar factor'. On the level of 
quadratic relations, such factors are certainly irrelevant and 
their discussion is beyond the scope of this text.\\[5pt]
\newlength{\parlength} \setlength{\parlength}{\textwidth}
\addtolength{\parlength}{-15em}
\begin{minipage}[b]{\parlength} \hspace*{1em} 
It is possible to give a more geometric interpretation 
of these formulas. To this end, let us consider an $N$-punctured
disc $D_N$. Its fundamental group $\pi_1(D_N)$ is freely generated 
by $N$ elements $l_\nu$. We can represent these generators by  
$N$ loops which surround the punctures as in Fig. 1. It was 
known to Artin already that the braid group $B_N$ acts on the  
fundamental group $\pi_1(D_N)$ and that this action is given by 
\ba
\s_\rho(l_\nu) \ = \ l_\nu \ \ & & \mbox{ for all } \ \ \nu \neq 
   \rho,\rho+1  \nn \\[2mm] 
\s_\rho(l_\rho) \ = \ l_{\rho+1} & , & \s_\rho(l_{\rho+1})\  = \    
        l_{\rho+1}^{-1}  l_\rho \ l_{\rho+1}. \nn
\ea 
\end{minipage} \ \ \raisebox{2ex}{
\begin{minipage}[b]{14em}\begin{center}
{\setlength{\unitlength}{0.00058333in}
\begin{picture}(3895,2514)(0,-10)
%\thicklines
\setlength{\maxovaldiam}{20pt}
\put(1800,1250){\oval(3650,2470)}
%\path(3883,12)(12,12)(12,2487)
%       (3883,2487)(3883,12)
\thicklines
\put(3079,689){\blacken\ellipse{84}{84}}
\put(3079,689){\ellipse{84}{84}}
\put(2445,1704){\blacken\ellipse{84}{84}}
\put(2445,1704){\ellipse{84}{84}}
\put(795,879){\blacken\ellipse{84}{84}}
\put(795,879){\ellipse{84}{84}}
\put(1956,2031){\makebox(0,0)[lb]{\smash{{{\SetFigFont{14}{16.8}{it}l}}}}}
\put(2083,1968){\makebox(0,0)[lb]{\smash{{{\SetFigFont{12}{14.4}{it}$\nu$}}}}}
\put(3364,1039){\makebox(0,0)[lb]{\smash{{{\SetFigFont{12}{14.4}{it}1}}}}}
\put(3237,1102){\makebox(0,0)[lb]{\smash{{{\SetFigFont{14}{16.8}{it}l}}}}}
\put(329,1345){\makebox(0,0)[lb]{\smash{{{\SetFigFont{14}{16.8}{it}l}}}}}
\put(456,1918){\makebox(0,0)[lb]{\smash{{{\SetFigFont{20}{24.0}{it}D}}}}}
\put(832,1779){\makebox(0,0)[lb]{\smash{{{\SetFigFont{12}{14.4}{it}N}}}}}
\put(456,1281){\makebox(0,0)[lb]{\smash{{{\SetFigFont{12}{14.4}{it}N}}}}}
\put(2873,1302){\blacken\ellipse{36}{36}}
\put(2873,1302){\ellipse{36}{36}}
\blacken\path(2779,1638)(2614,1445)(2694,1399)(2779,1638)
\path(2779,1638)(2614,1445)(2694,1399)(2779,1638)
\put(1890,736){\blacken\ellipse{116}{116}}
\put(1890,736){\ellipse{116}{116}}
\path(1908,722) (1848.585,760.842)
        (1791.041,797.985)
        (1735.309,833.452)
        (1681.332,867.268)
        (1629.052,899.458)
        (1578.412,930.043)
        (1529.354,959.050)
        (1481.819,986.501)
        (1435.751,1012.421)
        (1391.091,1036.833)
        (1347.781,1059.762)
        (1305.765,1081.231)
        (1264.983,1101.265)
        (1225.379,1119.886)
        (1186.895,1137.121)
        (1149.472,1152.991)
        (1077.582,1180.737)
        (1009.246,1203.315)
        (944.003,1220.918)
        (881.390,1233.736)
        (820.945,1241.962)
        (762.206,1245.786)
        (704.712,1245.402)
        (648.000,1241.000)

\path(648,1241) (600.453,1232.922)
        (549.373,1219.498)
        (497.389,1200.305)
        (447.132,1174.920)
        (401.235,1142.920)
        (362.326,1103.882)
        (333.037,1057.383)
        (316.000,1003.000)

\path(316,1003) (314.099,946.544)
        (326.302,893.476)
        (349.971,844.257)
        (382.470,799.348)
        (421.160,759.207)
        (463.403,724.295)
        (506.562,695.073)
        (548.000,672.000)

\path(548,672)  (599.521,648.693)
        (653.176,628.828)
        (709.466,612.422)
        (768.890,599.493)
        (831.945,590.059)
        (899.132,584.138)
        (970.949,581.746)
        (1008.750,581.880)
        (1047.895,582.903)
        (1088.448,584.817)
        (1130.469,587.624)
        (1174.023,591.328)
        (1219.171,595.930)
        (1265.976,601.432)
        (1314.499,607.836)
        (1364.804,615.145)
        (1416.953,623.361)
        (1471.007,632.486)
        (1527.031,642.522)
        (1585.085,653.471)
        (1645.232,665.337)
        (1707.534,678.120)
        (1772.055,691.824)
        (1838.856,706.450)
        (1908.000,722.000)

\blacken\path(757,1250)(993,1160)(1010,1252)(757,1250)
\path(757,1250)(993,1160)(1010,1252)(757,1250)
\put(1908,722){\blacken\ellipse{116}{116}}
\put(1908,722){\ellipse{116}{116}}
\path(1890,736) (1946.594,778.704)
        (2000.964,820.222)
        (2053.152,860.602)
        (2103.199,899.889)
        (2151.148,938.128)
        (2197.039,975.367)
        (2240.916,1011.651)
        (2282.819,1047.026)
        (2322.790,1081.538)
        (2360.871,1115.233)
        (2431.530,1180.357)
        (2495.131,1242.765)
        (2552.006,1302.826)
        (2602.491,1360.909)
        (2646.919,1417.380)
        (2685.624,1472.609)
        (2718.941,1526.964)
        (2747.203,1580.812)
        (2770.744,1634.522)
        (2789.898,1688.462)
        (2805.000,1743.000)

\path(2805,1743)        (2814.172,1790.408)
        (2819.556,1842.952)
        (2819.825,1898.308)
        (2813.653,1954.155)
        (2799.711,2008.170)
        (2776.673,2058.031)
        (2743.211,2101.415)
        (2698.000,2136.000)

\path(2698,2136)        (2645.444,2157.362)
        (2591.262,2164.132)
        (2536.793,2158.663)
        (2483.375,2143.304)
        (2432.348,2120.407)
        (2385.050,2092.323)
        (2342.821,2061.404)
        (2307.000,2030.000)

\path(2307,2030)        (2267.476,1989.677)
        (2230.459,1946.159)
        (2195.797,1898.973)
        (2163.336,1847.643)
        (2132.924,1791.697)
        (2104.409,1730.660)
        (2077.637,1664.058)
        (2052.457,1591.418)
        (2040.417,1552.684)
        (2028.716,1512.264)
        (2017.338,1470.096)
        (2006.262,1426.122)
        (1995.469,1380.284)
        (1984.941,1332.520)
        (1974.658,1282.773)
        (1964.601,1230.983)
        (1954.751,1177.090)
        (1945.089,1121.036)
        (1935.596,1062.761)
        (1926.254,1002.206)
        (1917.042,939.312)
        (1907.941,874.019)
        (1898.934,806.268)
        (1890.000,736.000)

\put(1681,1844){\blacken\ellipse{36}{36}}
\put(1681,1844){\ellipse{36}{36}}
\put(1357,1736){\blacken\ellipse{36}{36}}
\put(1357,1736){\ellipse{36}{36}}
\put(1087,1519){\blacken\ellipse{36}{36}}
\put(1087,1519){\ellipse{36}{36}}
\path(1913,732) (1978.332,704.399)
        (2041.520,678.128)
        (2102.626,653.175)
        (2161.710,629.525)
        (2218.833,607.167)
        (2274.056,586.086)
        (2327.439,566.271)
        (2379.043,547.708)
        (2428.930,530.384)
        (2477.160,514.286)
        (2523.793,499.401)
        (2568.891,485.716)
        (2612.514,473.219)
        (2654.723,461.895)
        (2695.579,451.732)
        (2735.142,442.717)
        (2773.474,434.838)
        (2810.636,428.080)
        (2881.689,417.879)
        (2948.788,412.010)
        (3012.420,410.368)
        (3073.071,412.851)
        (3131.227,419.353)
        (3187.375,429.771)
        (3242.000,444.000)

\path(3242,444) (3287.894,459.528)
        (3336.355,480.938)
        (3384.704,508.214)
        (3430.262,541.339)
        (3470.350,580.296)
        (3502.288,625.070)
        (3523.398,675.644)
        (3531.000,732.000)

\path(3531,732) (3523.108,788.623)
        (3501.755,839.286)
        (3469.639,884.023)
        (3429.456,922.867)
        (3383.904,955.854)
        (3335.679,983.016)
        (3287.479,1004.386)
        (3242.000,1020.000)

\path(3242,1020)        (3187.414,1033.937)
        (3131.304,1044.113)
        (3073.183,1050.421)
        (3012.564,1052.754)
        (2948.960,1051.008)
        (2881.884,1045.073)
        (2810.849,1034.846)
        (2773.694,1028.088)
        (2735.367,1020.217)
        (2695.807,1011.220)
        (2654.953,1001.083)
        (2612.743,989.792)
        (2569.118,977.334)
        (2524.015,963.697)
        (2477.375,948.866)
        (2429.136,932.829)
        (2379.238,915.572)
        (2327.619,897.081)
        (2274.218,877.345)
        (2218.976,856.348)
        (2161.831,834.078)
        (2102.721,810.521)
        (2041.587,785.665)
        (1978.367,759.496)
        (1913.000,732.000)

\blacken\path(3138,413)(2890,461)(2890,367)(3138,413)
\path(3138,413)(2890,461)(2890,367)(3138,413)
\put(1913,732){\blacken\ellipse{116}{116}}
\put(1913,732){\ellipse{116}{116}}
\end{picture} } \\[6.5mm]
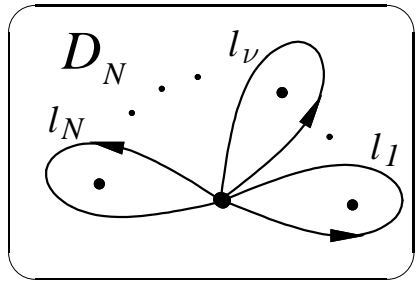
\parbox[b]{11.1em}{\footnotesize \sloppy {\bf Figure 1:} The fundamental
group $\pi_1(D_N)$ of an $N$-punctured disc $D_N$  is freely generated
by classes of the closed curves $l_\nu$.}
\end{center}
\end{minipage}} \\[2ex] 
\indent
Automorphisms of the fundamental group and of the multiloop-algebra
are denoted by the same letters. In any case, their action is
formally identical. We may push this even  further, if we assign
algebra valued matrices $\sM(p)$ to arbitrary elements $p \in 
\pi_1(D_N)$ by the obvious recursive definition 
$$\sM( p_1 p_2) \ :\sim \ \sM(p_1) \ \sM(p_2)
    \ \ \mbox{ for all } \ \ p_i \in  
     \pi_1(D_N) \ \ .$$
Notice that the objects $\sM(p)$ are products of the fundamental
matrices $\sM_\nu$ (up to scalar factors) and that they obey the
same exchange relations (\ref{M}) and transformation law under 
$\Phi_{\sf Ad}$ as the elements $\sM_\nu$. Using this
notation it is possible to lift the action of the braid group
on the fundamental group to an action on the multiloop-algebra
with the help of the formula
$$  \s (\sM(p))\  = \ \sM(\s(p))  \ \ . $$
Let us close this subsection by comparing Fig. 1 with the picture
we have drawn for $\cL_N$: while our pictorial presentations were
originally designed as a bookkeeping for the quantum symmetries 
of an algebra, they have turned out to encode all the classical 
symmetries at the same time.  
 
\subsection{Representations of the pure braid group $PB_N$}
The braid group $B_N$ contains an important subgroup $PB_N$ 
known as {\em pure braid group}. By definition, $PB_N$ is the 
kernel of the canonical homomorphism from the braid group into 
the symmetric group. One can show that $PB_N$ is generated by 
the elements $\eta_{\nu \mu}$ which are depicted in Fig. 2. 
In the context of multiloop-algebras, the pure braid group 
is singled out because its elements $\eta$ give rise to inner 
automorphisms of $\cL_N$. This means that for given $\eta \in 
PB_N$ there exists an element $\sv_\eta \in \cL_N$ such that 
$$    \eta( M) \ = \ \sv_\eta \ M \ \sv_\eta^{-1} \ \ 
      \mbox{ holds for all } \ \ M \in \cL_n\ \ . 
$$
Moreover, when $\eta$ is one of the standard generators 
$\eta_{\nu \mu} \in PB_N$, it is possible to give an  
explicit formula for the corresponding element $\sv_\eta
= \sv_{\nu \mu} $, 
$$ \sv_{\nu \mu} \ = \ V\/ [\, tr_q(\sM(l_\nu l_\mu))\,] \ \sim\ 
   V \/[\, tr_q(\sM_\nu \,\sM_\mu)\,]  \ \ \in \cL_N^\cG\ \ . $$
Here $V[.]$ is a universal function which does not depend 
on $\nu,\mu$. It is intimately related with the so-called 
{\em ribbon element} of $U_q(\sg)$ (see \cite{AlSc} for 
details).  

{\sloppy
\begin{prop} \cite{AlSc} {\em (Representations of $PB_N$)} There
exists an involution on $\cL^\cG_N$ such that the elements 
$\sv_\eta$ become unitary and the map $\eta \mapsto \sv_\eta 
\in \cL_N$ provides a unitary (projective) representation of 
the pure braid  group $PB_N$ by elements in the multiloop-algebra. 
\end{prop}}  \vspace*{-4ex}
\hspace*{-1em}
\begin{minipage}[b]{13em}
\begin{center} 
\setlength{\unitlength}{0.00033333in}
\begin{picture}(4578,4224)(0,-10)
\thicklines
\path(3428,782) (3427.764,829.001)
        (3427.685,863.565)
        (3428.000,907.000)

\path(3428,907) (3427.685,969.624)
        (3428.000,1032.000)

\path(3428,1032)        (3427.685,1095.000)
        (3428.000,1158.000)

\path(3428,1158)        (3429.016,1220.264)
        (3428.000,1283.000)

\path(3428,1283)        (3415.645,1346.336)
        (3396.000,1408.000)

\path(3396,1408)        (3357.361,1477.836)
        (3331.508,1518.370)
        (3302.636,1559.740)
        (3271.747,1599.768)
        (3239.843,1636.274)
        (3177.000,1690.000)

\path(3177,1690)        (3134.502,1713.743)
        (3088.716,1735.272)
        (3039.861,1754.728)
        (2988.157,1772.252)
        (2933.824,1787.982)
        (2877.083,1802.059)
        (2818.154,1814.623)
        (2757.255,1825.815)
        (2694.608,1835.775)
        (2630.432,1844.642)
        (2564.947,1852.557)
        (2498.374,1859.660)
        (2430.931,1866.091)
        (2362.840,1871.990)
        (2294.319,1877.498)
        (2225.590,1882.754)
        (2156.872,1887.899)
        (2088.384,1893.072)
        (2020.348,1898.415)
        (1952.983,1904.066)
        (1886.508,1910.167)
        (1821.144,1916.857)
        (1757.112,1924.276)
        (1694.630,1932.565)
        (1633.919,1941.864)
        (1575.198,1952.313)
        (1518.688,1964.052)
        (1464.609,1977.221)
        (1413.181,1991.960)
        (1364.623,2008.409)
        (1319.156,2026.709)
        (1277.000,2047.000)

\path(1277,2047)        (1209.685,2093.513)
        (1136.569,2160.409)
        (1070.418,2233.850)
        (1024.000,2300.000)

\path(1024,2300)        (1002.451,2361.665)
        (993.000,2426.000)

\path(993,2426) (1002.380,2490.238)
        (1024.000,2552.000)

\path(1024,2552)        (1070.470,2618.430)
        (1136.753,2691.887)
        %(1209.909,2758.651)
        %(1277.000,2805.000)

%\path(1277,2805)
%\path(1318.945,2825.619)
\path(1364.240,2844.432)
        (1412.662,2861.563)
        (1463.987,2877.136)
        (1517.994,2891.276)
        (1574.458,2904.107)
        (1633.157,2915.753)
        (1693.869,2926.337)
        (1756.369,2935.985)
        (1820.436,2944.821)
        (1885.846,2952.968)
        (1952.377,2960.551)
        (2019.805,2967.694)
        (2087.908,2974.521)
        (2156.463,2981.157)
        (2225.246,2987.725)
        (2294.036,2994.350)
        (2362.608,3001.156)
        (2430.740,3008.268)
        (2498.210,3015.809)
        (2564.793,3023.903)
        (2630.268,3032.676)
        (2694.412,3042.250)
        (2757.001,3052.750)
        (2817.812,3064.301)
        (2876.623,3077.027)
        (2933.211,3091.051)
        (2987.353,3106.498)
        (3038.826,3123.492)
        (3087.406,3142.158)
        (3132.872,3162.619)
        (3175.000,3185.000)

\path(3175,3185)        (3235.633,3232.752)
        (3299.176,3300.820)
        (3355.882,3374.228)
        (3396.000,3438.000)

\path(3396,3438)        (3416.054,3500.836)
        (3428.000,3565.000)

\path(3428,3565)        (3429.358,3627.951)
        (3428.000,3691.000)

\path(3428,3691)        (3428.000,3753.959)
        (3428.000,3817.000)

\path(3428,3817)        (3428.000,3880.526)
        (3428.000,3944.000)

\path(3428,3944)        (3428.000,3976.789)
        (3428.000,4039.000)

%\put(1277,2805){\blacken\ellipse{228}{228}}
%\put(1277,2805){\ellipse{228}{228}}
\path(12,783)(4566,783)(4566,625)
        (12,625)(12,783)
\path(12,783)(4566,783)(4566,625)
        (12,625)(12,783)
\path(4180,782)(4186,4039)
\path(2668,3185)(2668,4039)
\path(2676,782)(2676,1690)
\put(3212,72){\makebox(0,0)[lb]{\smash{{{\SetFigFont{14}{16.8}{rm}$\nu$}}}}}
\put(1137,72){\makebox(0,0)[lb]{\smash{{{\SetFigFont{14}{16.8}{rm}$\mu$}}}}}
\path(2676,1972)(2676,2944)
\put(2258,2426){\blacken\ellipse{60}{60}}
\put(2258,2426){\ellipse{60}{60}}
\put(1973,2426){\blacken\ellipse{60}{60}}
\put(1973,2426){\ellipse{60}{60}}
\put(1689,2426){\blacken\ellipse{60}{60}}
\put(1689,2426){\ellipse{60}{60}}
\path(12,4197)(4566,4197)(4566,4039)
        (12,4039)(12,4197)
\path(12,4197)(4566,4197)(4566,4039)
        (12,4039)(12,4197)
\path(513,782)(513,4040)
\path(1265,782)(1265,1941)
\path(1265,2160)(1265,4040)
\end{picture} \ \\[4mm] 
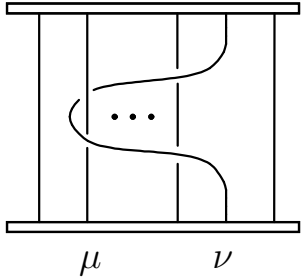
\parbox{11em}{\footnotesize \sloppy {\bf Figure 2:} Elements
$\eta_{\nu \mu}$ generate the pure braid group $PB_N$. They wrap
the $\nu^{th}$ strand once around the $\mu^{th}$. }
\end{center}
\end{minipage} \hspace*{2em}
\begin{minipage}[b]{16em}
\begin{center} 
\begin{picture}(160,160)(00,00) 
\setlength{\unitlength}{0.9pt}
{\setlength{\maxovaldiam}{20pt}\put(75,80){\oval(140,100)}}
\thicklines \linethickness{0.4mm} %\put(75,80){\oval(140,100)} 
{ \put(110,90){\oval(20,20)[t]} \put(50,90){\oval(20,20)[t]} 
\put(80,90){\oval(40,40)[b]}  \put(80,90){\oval(80,80)[b]}}
\put(50,90){\circle*{3}} \put(110,90){\circle*{3}} \put(80,90)
{\circle*{3}} \put(20,90){\circle*{3}} \put(140,90){\circle*{3}}
\put(25,90){$\dots$}  \put(63,90){$\dots$}  \put(85,90){$\dots$}
\put(123,90){$\dots$}{\Large { \put(52,82){$\mu$}  
\put(112,82){$\nu$}  \put(30,107){\huge $D_{\sss\hspace*{-3pt} N}$} 
\put(120,50){$l_{\nu \mu}$}}}
\end{picture} \\[-1ex] \hspace*{-2em}
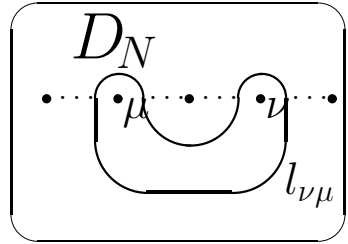
\parbox{14em}{\footnotesize \sloppy {\bf Figure 3:} Dehn twists
 along the curves $l_{\nu \mu}$ generate the mapping class group
 of $D_N$. They correspond to the elements $\eta_{\nu\mu}$ of the 
 pure braid group.} 
\end{center}
\end{minipage} \vspace*{3ex}

This simple result admits for more interesting generalizations
that we discuss in the next section. As a preparation we note
that the pure braid group may be reinterpreted as the mapping 
class group of the $N$-punctured disc. The mapping class group 
is the group of diffeomorphisms of a surface (in this case of 
$D_N$) modulo its identity component. It is generated by 
so-called Dehn twists which correspond to cutting the surface
along a given circle, rotating it by $2\pi$ and gluing back into  
the surface. When we perform this operation on a circle $l_\nu 
l_\mu$ in $D_N$, we obtain an element associated with the
generator $\eta_{\nu\mu}$ of the pure braid group (cf.
Fig. 2 and Fig. 3). Hence, we may rephrase the main result
this subsection by saying that our homomorphism $\eta \mapsto
\sv_\eta$ represents the mapping class group of the N-punctured
disc.
 
\section{Quantum Moduli Spaces and Mapping Class Groups}
\setcounter{equation}{0}
          
We are finally in a position to generalize the theory of the 
preceding section to surfaces of higher genus. This leads
us to the construction of quantum moduli spaces and to 
projective unitary representations of mapping class groups.  

\subsection{Theory on a torus} It is instructive to study the 
torus $T$ first. It has two nontrivial cycles, the $a$- and the
$b$-cycle, which generate the fundamental group $\pi_1(T)$. 
For simplicity, let us remove a disc $D$ from the torus\ $T$ and look \\[2mm]
\begin{minipage}[b]{\parlength} 
at the group $\pi_1(T\setminus D)$ which is freely generated 
by  $a,b$. The mapping class group of $T\setminus D$ is 
generated by two Dehn twists $\a,\b$ along the 
$a,b$-cycle. They act on the fundamental group
according to 
\ba 
  \a(a) \ =\  a \ \ \ \ & \ , \ & \a (b)\  = \ b\, a \nn\ \ , \\[2mm]
  \b(a) \ =\  b^{-1} \, a  & , & \b(b) \ = \ b \nn\ \ .  
\ea
This action obviously respects the standard commutator
relation $[a,b]:= b\, a^{-1}b^{-1}a=e$ and hence descend 
to the fundamental group $\pi_1(T)$ of the torus.  
\end{minipage} \raisebox{7ex}{
\begin{minipage}[b]{14em}
\begin{center} 
\setlength{\unitlength}{0.00053333in}
\begin{picture}(3439,3454)(0,-10)
\linethickness{2pt}
\path(1720,3427)(3427,1719)(1720,12)
        (12,1719)(1720,3427)(1720,3427)
%\path(1044,2715)(1152,2857)(1008,2751)
%\path(2431,687)(2289,580)(2397,722)
%\path(2397,2715)(2289,2857)(2431,2751)
%\path(1044,722)(1152,580)(1008,687)
\thicklines
\path(2573,866)(866,2573)
\blacken\path(2077,2003)(2218,2216)(2006,2074)(2077,2003)
\path(2077,2003)(2218,2216)(2006,2074)(2077,2003)
\blacken\path(1436,2074)(1223,2216)(1365,2003)(1436,2074)
\path(1436,2074)(1223,2216)(1365,2003)(1436,2074)
\put(1288,1026){\makebox(0,0)[lb]{\smash{{{\SetFigFont{14}{16.8}{it}a}}}}}
\put(2313,1201){\makebox(0,0)[lb]{\smash{{{\SetFigFont{14}{16.8}{it}b}}}}}
\put(1313,2526){\makebox(0,0)[lb]{\smash{{{\SetFigFont{20}{24.0}{it}T}}}}}
\path(866,866)(2573,2573)
\put(1720,1719){\blacken\ellipse{100}{100}}
\put(1720,1719){\ellipse{100}{100}}
\end{picture} 
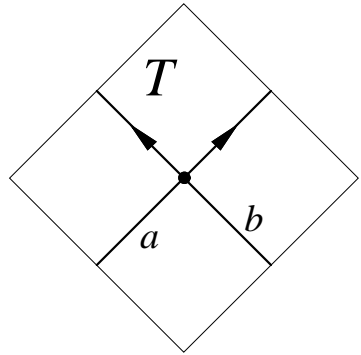
\parbox[t]{11em}{\footnotesize \sloppy {\bf Figure 4:} The $a-$ and 
$b$-cycle generate the fundamental group of $T$ and $T \setminus 
D$. Opposite edges are identified.}  
\end{center}
\end{minipage}} \vspace*{1ex}   

Next we have to understand what replaces the multiloop-algebra 
in the genus 1 example we are dealing with now. The basic idea
is to assign two copies of the loop-algebra $\cL$ to the two 
generators $a,b$ of the fundamental group. They will be denoted 
by $\sA = \sM (a)$ and $\sB = \sM (b)$. We want to combine them 
so that the resulting algebra $\cH$ preserves the quantum symmetry 
of the loop-algebras and admits an action of the mapping class 
group. The second requirement implies that
\ba      \a (\sA) \ = \ \sA \ \ \ \ & \ \ ,\ \    & \a (\sB) \
         \sim \ \sB\  \sA \ \ ,  \\[2mm]
        \b(\sA) \ \sim\  \sB^{-1} \sA  & , & \b(\sB) \ = \ \sB  
\ea 
extend to automorphisms of $\cH$. Observe that the formulas for 
the action of $\a,\b$ on $\cH$ are obtained by lifting the action 
of the mapping class group on the fundamental group. 
 
With these to basic requirements one is lead to an algebra $\cH$ 
which is generated by matrix elements of $\sA,\sB$ subject to 
the defining relations of $\cL$ and the exchange relation
\ba 
    \rR^{-1}\ \Ae{1}{} \ \ \rR \ \Be{2}{} & = & \Be{2}{}\ \  
      \rR' \ \Ae{1}{} \ \ \rR
     \\[2mm]
    \cH \  \mbox{ has QS: } \ \  \Phi_{\sf Ad}:\  \sA
    \ \mapsto & & \hspace*{-1.5em} \sT^{-1} \sA \ \sT 
    \ \  , \ \  \Phi_{\sf Ad}: \ B \ \mapsto \ \sT^{-1} \sB \ \sT 
    \ \ . \nn 
\ea 
We call $\cH$ the {\em handle-algebra}. A more detailed investigation 
reveals that elements of the mapping class group act as inner 
automorphisms on the handle-algebra. The action of $\a, \b$ is 
implemented by conjugation with the elements 
$$    
     \sv_\a \ = \ V\/ [\, \tr _q(\sA)\,] \ \ \ \mbox{and} \ \ \ 
      \sv_\b \ = \ V\/ [\, \tr _q(\sB)\,]\ \ . 
$$                         
Here $V$ is the same universal function we encountered in
Subsection 3.2. 
It can be shown that the map $\a \mapsto \sv_\a$,
$ \b \mapsto \sv_\b$ gives rise to a projective unitary representation
of the mapping class group of $T \setminus D$. One may consistently
impose the constraint $ \sB \sA^{-1} \sB^{-1} \sA \sim {\bf 1}$ to
obtain a theory for the torus $T$.  

\subsection{Theory for arbitrary genus} It should be clear by now, 
how to proceed in the general case of an arbitrary Riemann surface
$\Sigma_{g,N}$ of genus $g$ and with $N$ punctures. The fundamental 
group of $\Sigma_{g,N} \setminus D$ is freely generated by $2g+N$
cycles $l_\nu, a_i, b_i$. Within $\pi_1(\Sigma_{g,N})$, these generators
obey one relation 
$$ r \ = \ [\,a_g\, ,\,b_g\,]\  \cdots\  [\,a_1
\, ,\, b_1\,]\  l_N \ \cdots \ l_1\  = \ e \ \ . $$ 
As in the two examples we have studied, the fundamental group 
$\pi_1(\Sigma_{g,N} \setminus D)$ carries an action of the mapping class 
group of $\Sigma_{g,N} \setminus D$. Explicit formulas can be worked out.

Let us define an algebra $\cL_{g,N}$ to be the braided tensor product 
of $N$ copies of the loop-algebra $\cL$ and $g$ copies of the 
handle-algebra $\cH$ so that the generating matrices 
$\sM_\nu$,  $\sA_i$ and $\sB_i$ ($i = 1, \dots, g$) obey the following 
additional exchange relations      
\ba 
 \rR^{-1} \M{1}_\nu\  \rR \ \Ae{2}_i  & = &  
     \Ae{2}_i \ \rR^{-1} \M{1}_\nu \  \rR\ \ \  
      \ \ \mbox{ for all } \ \ i,\nu  \\[2mm]
   \rR^{-1} \Ae{1}_i\  \rR \ \Ae{2}_j  & = &  
     \Ae{2}_j \ \rR^{-1} \Ae{1}_i \  \rR\ \ \  
      \ \  \ \ \ i < j  
\ea
and similarly when $\sA_i$ and/or $\sA_j$ are replaced by $\sB_i, 
\sB_j$.   
By construction, $\cL_{g,N}$ contains $2g +N$ copies of $\cL$ which 
we think of as being assigned to the $2g+N$ fundamental cycles of 
the surface. This allows to lift the action of the mapping class 
group from the fundamental group to the algebra $\cL_{g,N}$. All 
the automorphisms of $\cL_{g,N}$ which arise in this way are inner 
and we end up with a projective (unitary) representation of the 
mapping class group by elements in the algebra $\cL^\cG_{g,N}$
of elements transforming trivially under $\Phi_{\sf Ad}$.   
   
\begin{prop} \cite{AlSc} {\em (Representations of the mapping
class group)} Let $p$ denote an arbitrary cycle on $\Sigma_{g,N}
\setminus D$ and $\gamma_p$ the associated Dehn twist in the
mapping class group. Then there exists an involution on 
$\cL_{g,N}^\cG$ such that
$$   \gamma_p \ \mapsto \ \sv (\gamma_p) \ : = \
     V\/ [\, \tr_q ( M(p) )\,] \ \ \in \ \cL^\cG _{g,N} \ \subset 
     \ \cL_{g,N} 
$$
extends to a projective unitary representation of the mapping
class group of $\Sigma_{g,N} \setminus D$. Here $M(p)$ is
defined recursively as in Subsection 3.2. The map $\gamma_p
\mapsto \sv (\gamma_p)$ descends to the mapping class group 
of $\Sigma_{g,N}$, if $\cL^\cG_{g,N}$ is factored by the 
relation $\sM(r) \sim {\bf 1}$ .      
\end{prop}

Representations of the algebras $\cL_{g,N}$ and $\cL_{g,N}^\cG$ 
have been constructed in \cite{AlSc}. They give rise to projective 
unitary actions of mapping class groups on finite dimensional 
Hilbert spaces. The latter are equivalent to the representations 
obtained by Reshetikhin and Turaev \cite{ReTu}. 
 
We are finally able to explain how the quantum moduli spaces 
$\cA_{\sf CS}$ arise within the described framework. The idea 
it to descend from the algebra $\cL_{g,N}$ generated by the 
quantum connections in two steps. First one restricts to the 
subspace of elements on which the quantum gauge symmetry 
$\Phi_{\sf Ad}$ acts trivially. This simply gives our algebra 
$\cL^\cG_{g,N}$. Then, in a second step, the flatness condition 
$\sM (r) \sim {\bf 1}$ is imposed, i.e. 
$$     \cA_{\sf CS} \ := \  \cL_{g,N}^\cG \ / \langle \, \sM (r)
       \, \sim \, {\bf 1}\, \rangle \ \ . $$
It was shown in \cite{AlSc} that this algebra $\cA_{\sf CS}$ 
possesses irreducible representations on the spaces of 
conformal blocks of the corresponding WZNW-model. This is 
in perfect agreement with the results of \cite{geo}.

\end{document}